\input amstex
\documentstyle{amsppt}
\magnification 1200
\NoRunningHeads
\NoBlackBoxes
\document

\def\Ua{U_q(\tilde\g)}
\def\U2{{\Ua}_2}
\def\g{\frak g}

\def\Z{\Bbb Z}
\def\C{\Bbb C}

\def\l{\lambda}

\def\<{\langle}
\def\>{\rangle}
\def\o{\otimes}

\def\End{\text{End}}

\topmatter
\title On set-theoretical solutions of the quantum Yang-Baxter equation
\endtitle

\author Pavel Etingof, Travis Schedler, and Alexandre Soloviev
\endauthor

\address
\newline
P.~E. : Department of Mathematics, Harvard University, Cambridge, MA
02138
\newline
T.~S. : 1500 W.Sullivan Rd., Aurora, IL 60506
\newline
A.~S. : Department of Mathematics, MIT, Cambridge, MA 02139
\endaddress

\email
\newline
P.~E. : etingof\@ math.harvard.edu
\newline
T.~S. : trasched\@ imsa.edu
\newline
V.~R. : sashas\@ math.mit.edu
\endemail

\endtopmatter

In the paper  \cite{Dr}
 V.Drinfeld formulated a number of problems in quantum group theory.
In particular, he suggested to consider ``set-theoretical'' solutions 
of the quantum Yang-Baxter equation, i.e. solutions given by a permutation 
$R$ 
of the set $X\times X$, where $X$ is a fixed 
finite set. In this note we study such 
solutions, which satisfy the unitarity 
and the crossing symmetry conditions -- natural conditions arising 
in physical applications. More specifically, we consider 
``linear'' solutions: the set $X$ is an abelian group, 
and the map $R$ is an automorphism of $X\times X$. 
We show that in this case, 
solutions are in 1-1 correspondence with pairs $a,b\in \End X$, 
such that $b$ is invertible and $bab^{-1}=\frac{a}{a+1}$. 
Later we consider ``affine'' solutions ($R$ is an automorphism 
of $X\times X$ as a principal homogeneous space), and show that they have a 
similar classification. The fact that these classifications are so nice 
leads us to think that there should be some interesting structure 
hidden behind this problem. 

Let $X$ be a finite set, and $R:X\times X\to X\times X$ be a bijection. 
Recall \cite{Dr} that the quantum Yang-Baxter equation 
is the equation 
$$
R^{12}R^{13}R^{23}=R^{23}R^{13}R^{12},\tag 1
$$
and the unitarity condition is
$$
R^{21}R=1. \tag 2
$$

We will look for solutions of (1) and (2) of the following form:
$X$ is an abelian group, and $R$ is an automorphism of $X\times X$
(``linear'' solutions). This was inspired by the paper \cite{Hi}. 
In this case $R$ has the form 
$$
R(x,y)=(cx+dy,ax+by),\ a,b,c,d\in \End X. \tag 3
$$
It is easy to check that for $R$ of the form (3) equation (1) is equivalent 
to the equations \cite{Hi}
$$
\gather
a(1-a)=bac,\ d(1-d)=cdb,\ ab=ba(1-d),\ ca=(1-d)ac,\ dc=cd(1-a),\\ 
bd=(1-a)db,\  cb-bc=ada-dad.\tag 4
\endgather
$$
It is also easy to see that equation (2) is equivalent to the equations
$$
a^2+bc=1,\ cb+d^2=1,\ ab+bd=0,ca+dc=0.\tag 5
$$

Now consider the crossing symmetry condition. If $R$ satisfies (1) and (2), 
the crossing symmetry condition is
$$
(((R^{-1})^t)^{-1})^t=R,\tag 6
$$
where $()^t$ denotes transposing in the second component
of the product (here $R$ is regarded as a matrix of 0-s and 1-s). 

Using (2), we can rewrite condition (6) as 
$$
(R^{21})^tR^t=1.\tag 7
$$
Since $R=\sum_{x,y}E_{cx+dy,x}\o E_{ax+by,y}$, where $E_{pq}$ is 
the elementary matrix, condition (7) can be written in the form
$$
\sum_{y'=cx+dy,y=cx'+dy'}E_{ax'+by',x}\o E_{x',ax+by}=1.\tag 8
$$
But $1=\sum_{p,q}E_{pp}\o E_{qq}$. This implies that equation (6), 
given (1),(2), is equivalent to the following:

1. For any $x,x'\in X$ there exist unique $y,y'\in X$ such that 
$y'=cx+dy,y=cx'+dy'$, and

2. These $y,y'$ satisfy the equations $x'=ax+by,x=ax'+by'$. 

The first condition is equivalent to the condition that the matrix 
$\left(\matrix 1&-d\\ -d&1\endmatrix\right)$ is invertible, i.e. 
that $1-d^2$ is invertible. But by (5) $1-d^2=cb$, so $c,b$ are invertible. 

\proclaim{Proposition 1} If $b,c$ are invertible, equations (4),(5) are 
equivalent to the equations 
$$
bab^{-1}=\frac{a}{a+1},c=b^{-1}(1-a^2),d=\frac{a}{a-1}.\tag 9
$$
\endproclaim

\demo{Proof} The second equation of (9) follows directly from (5). 
Also, (5) implies 
$$
a=-bdb^{-1}.\tag 10
$$
Therefore, multiplying the equation $bd=(1-a)db$ (which is in (4)) by $b^{-1}$ 
on the right, we get 
$$
-a=(1-a)d.\tag 11
$$
Since $b,c$ are invertible, so is $bc=1-a^2$, so 
$1-a$ is invertible. Thus, (11) implies the third equation of (9). 
Now the first equation of (9) follows from (10). 

Conversely, substituting (9) into (4),(5), it is easy to show 
by a direct calculation that they are identically satisfied. 
$\square$\enddemo

\proclaim{Corollary 2} A map $R$ of the form (3) is a solution to (1),(2),(6) 
if and only if $b,c$ are invertible, and (9) are satisfied. 
Thus, such solutions are in 1-1 correspondence with pairs (a,b) such 
that $bab^{-1}=\frac{a}{a+1}$. 
\endproclaim

\demo{Proof} It remains to show that if $b,c$ are invertible and (9) 
are satisfied then condition 2 holds. This is checked by a direct calculation. 
\enddemo

Now consider the case when $R$ has the form 
$$
R(x,y)=(cx+dy+t,ax+by+z), t,z\in X.\tag 12
$$
(``affine'' solutions, cf \cite{Hi}).
In this case, it is clear that the equations on $a,b,c,d$ are the same as 
before. The only equation for $z,t$ is obtained from (1) and has the form
$t=-b^{-1}(1+a)z$. Thus, we get

\proclaim{Proposition 3} 
A map $R$ of the form (12) is a solution to (1),(2),(6) 
if and only if $b,c$ are invertible, (9) are satisfied, and
$t=-b^{-1}(1+a)z$. 
Thus, such solutions are in 1-1 correspondence with triples (a,b,z) such 
that $bab^{-1}=\frac{a}{a+1}$. 
\endproclaim

Now consider examples of solutions of the equation 
$$
bab^{-1}=\frac{a}{a+1}.\tag 13
$$ 

{\bf Example 1.} \cite{Hi} Let $X=\Z/n \Z$. Then $\End X=\Z/n\Z$, 
which is commutative, so equation (13) reads $a=\frac{a}{a+1}$, which is 
equivalent to $a^2=0$ (and $b$ is any invertible element). 

{\bf Example 2.} Let $X=V^N$, where $V$ is an abelian group.
Then the algebra $Mat_N(\Z)$ of integer matrices is mapped into 
$\End X$. Thus, it is enough for us to construct a solution of 
$bab^{-1}=\frac{a}{a+1}$ in $Mat_N(\Z)$, such that $b\in GL_N(\Z)$. 
 
Let $a_{ij}=\delta_{i+1,j}$, and $b_{ij}=\left(\matrix j\\ i\endmatrix\right)$. 
Then $a,b$ satisfy (13). Indeed, this equation can be 
rewritten as $ab=ba+aba$, which at the level of matrix elements
reduces to the well-known identity
for binomial coefficients:
$$
\left(\matrix j\\ i+1\endmatrix\right) =
\left(\matrix j-1\\ i\endmatrix\right) +
\left(\matrix j-1\\ i+1\endmatrix\right) \tag 14
$$ 
We will use the following notation for this solution: 
$a=J_N$, $b=B_N$. 

In fact, all solutions of (13) in $Mat_N(\Z)$ can be obtained 
from this one. Indeed, we have 

\proclaim{Lemma 4} Let $a,b$ be a solution of (13) in $Mat_N(\C)$. 
Then $a$ is nilpotent. 
\endproclaim

\demo{Proof} It follows from (13) that if $\l$ is an eigenvalue 
of $a$ then so is $\frac{\l}{\l+1}$. Therefore, if $\l\ne 0$, 
we get that $a$ has infinitely many distinct eigenvalues. 
This is impossible, so $\l=0$.
$\square$\enddemo

Thus, $a$ is nilpotent. Then,
by Jordan's theorem, $a$ can be reduced, over $\Z$, to Jordan normal form:
$a=J_{N_1}\oplus...\oplus J_{N_K}$, where 
$J_{N_l}\in Mat_{N_l}(\Z)$ are given by $(J_{N_l})_{ij}=\delta_{i+1,j}$. 
If $a$ is of this form, then $b=b_0A$, where $A$ commutes with $a$, and 
$b_0=B_{N_1}\oplus...\oplus B_{N_K}$.
Thus we have proved 

\proclaim{Proposition 5} Any solution of (13) in $Mat_N(\Z)$ with $b\in 
GL_N(\Z)$ is conjugate under $GL_N(\Z)$ to a solution of the form
$a=J_{N_1}\oplus...\oplus J_{N_K}$,
$b=(B_{N_1}\oplus...\oplus B_{N_K})A$, where $[A,a]=0$. 
\endproclaim

\proclaim{Proposition 6} If $V=\Z/p\Z$, where $p$ is a prime, 
and $N<p$, then any solution of (13) in $\End X$ is 
conjugate to a solution of the form 
given in Proposition 5. 
\endproclaim

\demo{Proof} Let $\l$ be an eigenvalue of $a$ over 
the algebraic closure $\overline{\Z/p\Z}$. 
Then $\frac{\l}{\l+1}$ is also an eigenvalue. Therefore, 
if $\l\ne 0$, $a$ has to have at least $p$ distinct eigenvalues. 
The rest of the proof is the same as for Proposition 5. 
$\square$\enddemo

However, if $N\ge p$, other solutions are possible. 

{\bf Example: } $p=N=2$, $a=\left(\matrix 1&1\\ 1&0\endmatrix\right)$, 
$b=\left(\matrix 0&1\\ 1&0\endmatrix\right)$.

\enddocument